\documentclass [a4paper,10pt]{article}
\usepackage{amssymb,amsmath,amscd,url}
\newtheorem{defi}{Definition}[section]
\newtheorem{theo}[defi]{Theorem}
\newtheorem{rema}[defi]{Remark}
\DeclareMathOperator{\trace}{trace}
\DeclareMathOperator{\divergence}{div}

\title{The Fujita phenomenon in exterior domains under the Robin boundary conditions} 
\author{Jean-Fran\c cois Rault}
\date{}

\begin{document}

\maketitle

\begin{center}
LMPA Joseph Liouville, ULCO, FR 2956 CNRS, \\
Universit\'e Lille Nord de France \\
50 rue F. Buisson, B.P. 699, F-62228 Calais Cedex (France)\\
\url{jfrault@lmpa.univ-littoral.fr}\\
\end{center}
\bigskip

\begin{abstract}
\noindent The Fujita phenomenon for nonlinear parabolic problems $\partial_t u=\Delta u +u^p$ in an exterior domain of $\mathbb{R}^N$ under the Robin boundary conditions is investigated in the superlinear case.  As in the case of Dirichlet boundary conditions (see Trans. Amer. Math. Soc 316 (1989), 595-622 and Isr. J. Math. 98 (1997), 141-156), it turns out that there exists a critical exponent $p=1+2/N$ such that blow-up of positive solutions always occurs for subcritical exponents, whereas in the supercritical case global existence can occur for small non-negative initial data.\\

 
\end{abstract}
\bigskip
\section{Introduction}

Let $\Omega$ be an exterior domain of $\mathbb{R}^N$ , that is to say a connected open set $\Omega$ such that $\overline{\Omega}^c$ is a bounded domain when $N\geq 2$, and in dimension one, $\Omega$ is the complement of a real closed interval. We always suppose that the boundary $\partial \Omega $ is of class $\mathcal{C}^2$. The outer normal unit vector field is denoted by $\nu: \partial \Omega \rightarrow \mathbb{R}^N $ and the outer normal derivative by $\partial _\nu$. Let $p$ be a real number with $p>1$, $\alpha$ a non-negative continuous function on $\partial \Omega \times \mathbb{R}^+$ and $\varphi$ a continuous function in $\overline{\Omega}$. Consider the following nonlinear parabolic problem
\begin{equation}\label{Ext_Robin}
\left\{
    \begin{array}{ll}
        \partial _t u = \Delta u + u^p& \textrm{ in }\overline{\Omega} \times (0,+\infty), \\
       \partial _\nu u+\alpha u= 0 & \textrm{ on } \partial \Omega \times (0,+\infty), \\
        u(\cdot , 0) = \varphi & \textrm{ in } \overline{\Omega} .
    \end{array}
\right.
\end{equation}
In this paper, we give a positive answer to Levine \& Zhang's question \cite{LZ}: the Fujita phenomenon, well-known in the case of $\Omega=\mathbb{R}^N$ (see Ref. \cite{Fujita}), remains true for the Robin boundary conditions. The case limiting $\alpha \equiv 0$ and $\alpha = +\infty$ were proved by Levine \& Zhang in \cite{LZ} and by Bandle \& Levine in \cite{BL1}, respectively. The real number $1+\frac{2}{N}$ is still the critical exponent, and we prove the blowing-up of all positive solutions of Problem $(\ref{Ext_Robin})$ for subcritical exponents $p$, whereas in the supercritical case, we show the existence of global positive solutions of Problem $(\ref{Ext_Robin})$ for sufficiently small initial data. In the last section, we study the case of a general second order elliptic operator replacing the Laplacian. We also consider a non-linearity including a time and a space dependence. Throughout, we shall assume that $\alpha$ is non-negative
\begin{equation}\label{H0_Robin}
\alpha \geq 0 \textrm{ on } \partial \Omega \times \mathbb{R}^+ ,
\end{equation}
and, in order to deal with classical solutions, we need some regularity on $\alpha$
\begin{equation}\label{H1_Robin}
\alpha \in \mathcal{C}(\partial \Omega \times \mathbb{R}^+).
\end{equation}
To construct solutions with the truncation procedure (see Section \ref{Prelem}), we suppose 
\begin{equation}\label{i-data}
\varphi \in \mathcal{C}(\overline{\Omega}), \ 0<\parallel \varphi \parallel_\infty < \infty, \ \varphi \geq 0,\ \lim_{\parallel x \parallel_2 \rightarrow\infty}\varphi (x)=0 .
\end{equation}
In the case $\Omega = \mathbb{R}^N$, the boundary conditions are dropped, and the result is well-known by  the classical paper of Fujita \cite{Fujita}. Thus we suppose $\Omega \not= \mathbb{R}^N$.

\section{Preliminaries}\label{Prelem}

First, we give the definition of positive solution which is understood along this paper.
\begin{defi}
A positive solution of Problem $(\ref{Ext_Robin})$ is a positive function $u: (x,t) \mapsto u(x,t)$ of class $\mathcal{C}(\overline{\Omega} \times [0,T)) \cap \mathcal{C}^{2,1}(\overline{\Omega} \times (0,T))$, satisfying 
\begin{equation*}
\left\{
    \begin{array}{ll}
        \partial _t u = \Delta u + u^p& \textrm{ in }\overline{\Omega} \times (0,+\infty), \\
       \partial _\nu u+\alpha u= 0 & \textrm{ on } \partial \Omega \times (0,+\infty), \\
        u(\cdot , 0) = \varphi & \textrm{ in } \overline{\Omega} ,
    \end{array}
\right.
\end{equation*}
where $\alpha$ and $\varphi$ are given with $(\ref{H0_Robin})$, $(\ref{H1_Robin})$ and $(\ref{i-data})$. The time $T=T(\alpha, \varphi) \in (0,+\infty]$ denotes the maximal existence time of the solution $u$. If $T=+\infty$, the solution is called global.
\end{defi}
From \cite{BL1}, if $T<+\infty$, $u$ blows up in finite time, that is to say:
\begin{displaymath}
\lim_{t\nearrow T} \sup_{x\in \overline{\Omega}} u(x,t)= +\infty.
\end{displaymath}
Then, let us recall a standard procedure to construct solutions of Problem $(\ref{Ext_Robin})$ in outer domains for uniformly bounded and continuous initial data $\varphi$. For more details, we refer to \cite{BBR1} , \cite{JFR} and references therein. Let $( D_n)_{n \in \mathbb{N}}$ be a sequence of nested bounded domains such that
\begin{displaymath}
\overline{\Omega}^c \subseteq D_0 \subseteq D_1 \subseteq \dots \subseteq \bigcup_{n \in \mathbb{N}} D_ n = \mathbb{R}^N .
\end{displaymath}
Let $u_n$ be the solution of 
\begin{equation}\label{truncP}
\left\{
    \begin{array}{ll}
        \partial _t u = \Delta u + u^p& \textrm{ in }\overline{\Omega} \cap D_n \times (0,+\infty), \\
       \partial _\nu u+\alpha u= 0 & \textrm{ on } \partial \Omega \times (0,+\infty), \\
       u= 0 & \textrm{ on } \partial D_n \times (0,+\infty), \\
        u(\cdot , 0) = \varphi_n & \textrm{ in } \overline{\Omega} \cap D_n ,
    \end{array}
\right.
\end{equation}
where $(\varphi_n)_{n \in \mathbb{N}}$ denotes a sequence of functions in $\mathcal{C}_0(\overline{\Omega} \cap D_n)$ such that
\begin{displaymath}
0 \leq \varphi_n \leq \varphi \textrm{ in }  \overline{\Omega} \cap D_n
\end{displaymath}
and $\varphi_n \to \varphi$ uniformly in any compact of $\overline{\Omega} \cap D_n$ as $n \to +\infty$. Let $z$ denote the solution of the ODE 
\begin{equation*}
\left\{
    \begin{array}{ll}
        \dot{z} = z^p, \\
        z(0)= \parallel \varphi \parallel_\infty,
    \end{array}
\right.
\end{equation*}
with maximal existence time $S=\frac{1}{(p-1)\parallel \varphi \parallel_\infty^{p-1}}$. By the comparison principle (see \cite{BDC1}), we have
\begin{displaymath}
0 \leq u_n(x,t) \leq u_{n+1}(x,t) \leq z(x,t) \textrm{ in }  \overline{\Omega} \cap D_n \times [0,S].
\end{displaymath}
Standard arguments based on a priori estimates for the heat equation imply  $u_n \rightarrow u$ in the sense of $\mathcal{C}^{2,1}_{loc}(\overline{\Omega} \times (0,S))$ as $n\rightarrow +\infty$, where $u$ is a positive solution of Problem $(\ref{Ext_Robin})$. Moreover, since $u_n$ vanishes on $\partial D_n$ for each $n\in \mathbb{N}^*$, the solution $u$ vanishes at infinity: 
\begin{displaymath}
\lim_{\parallel x \parallel_2 \rightarrow \infty} u(x,t) = 0 \ , \forall \ t \in (0,T)  .
\end{displaymath}

\section{Blow up case}

In this section, we compare the solution of Problem $(\ref{Ext_Robin})$ with an appropriate Dirichlet solution. We prove the following theorem:

\begin{theo}\label{explode} Suppose that conditions $(\ref{H0_Robin})$, $(\ref{H1_Robin})$ and $(\ref{i-data})$ are fullfiled. Then all non-trivial positive solutions of Problem $(\ref{Ext_Robin})$ blow up in finite time for $p \in (1, 1+ 2/N)$. Moreover, if $N \geq 3$, blow up also occurs for $p=1 + 2/N$.
\end{theo}
\emph{Proof: } Ab absurdo, suppose that there exists $\alpha$ and a non-trivial $\varphi$ satisfying the hypotheses above, and such that the solution $u$ of Problem $(\ref{Ext_Robin})$  with these parameters is global. Then, consider $u_n$ the solution of the truncated Problem $(\ref{truncP})$. By the comparison principle from \cite{BDC1}, we obtain
\begin{displaymath}
0 \leq u_n(x,t) \leq u(x,t) \textrm{ in }  \overline{\Omega} \cap D_n \textrm{ for } t >0.
\end{displaymath}
Thus, $u_n$ can not blow up in finite time, and $u_n$ must be global. Next, define $v_n$ the solution of the following problem
\begin{equation*}
\left\{
    \begin{array}{ll}
        \partial _t v_n = \Delta v_n + v_n^p& \textrm{ in }\overline{\Omega} \cap D_n \times (0,+\infty), \\
       v_n= 0 & \textrm{ on } \partial \Omega \times (0,+\infty), \\
       v_n= 0 & \textrm{ on } \partial D_n \times (0,+\infty), \\
       v_n(\cdot , 0) = \varphi_n & \textrm{ in } \overline{\Omega} \cap D_n .
    \end{array}
\right.
\end{equation*}
Again, the comparison principle from \cite{BDC1} implies $0 \leq v_n(x,t) \leq u_n(x,t)$ in $\overline{\Omega} \cap D_n$ for  $t >0$. Then, we consider $v$ the solution of the Dirichlet problem
\begin{equation*}
\left\{
    \begin{array}{ll}
        \partial _t v = \Delta v + v^p& \textrm{ in }\overline{\Omega} \times (0,+\infty), \\
       v= 0 & \textrm{ on } \partial \Omega \times (0,+\infty), \\
       v(\cdot , 0) = \varphi & \textrm{ in } \overline{\Omega}  ,
    \end{array}
\right.
\end{equation*}
obtained as the limit of the $v_n$ by the truncation procedure described in Section \ref{Prelem}. Thus, $v \leq u$ in $\overline{\Omega} \times (0,+\infty)$ and $v$ is a global positive solution. A contradiction with Bandle \& Levine results \cite{BL1} (see \cite{BL2} for the one-dimensional case). If $N \geq 3$ and $p=1 + 2/N$, the contradiction holds with Mochizuki \& Suzuki's results \cite{MS} and \cite{Suzuki}. Hence, our solution $u$ must blow up in finite time. 
\\
\mbox{}\nolinebreak\hfill\rule{2mm}{2mm}\par\medbreak

\section{Global existence case}

Now, we consider supercritical exponents:
\begin{displaymath}
p > 1 + \frac{2}{N}.
\end{displaymath}
We look for a global positive super-solution of Problem $(\ref{Ext_Robin})$, we mean a function $U$ satisfying
\begin{equation*}
\left\{
    \begin{array}{ll}
        \partial _t u \geq \Delta u + u^p& \textrm{ in }\overline{\Omega} \times (0,+\infty), \\
       \partial_\nu u+\alpha u \geq 0 & \textrm{ on } \partial \Omega \times (0,+\infty), \\
        u(\cdot , 0) \geq \varphi & \textrm{ in } \overline{\Omega} .
    \end{array}
\right.
\end{equation*}
With this global super-solution and using the comparison principle, we construct the sequence $(u_n)_{n \in \mathbb{N}}$ of global positive solutions of Problems $(\ref{truncP})$. Thus, using the truncation procedure of Section \ref{Prelem}, we construct a global positive solution of Problem $(\ref{Ext_Robin})$. We use two different super-solutions, and we obtain two results on the global existence with some restrictions on the dimension $N$ or on the coefficient $\alpha$. First, we only suppose that the dimension 
\begin{displaymath}
N \geq 3.
\end{displaymath}

\begin{theo}
Under hypotheses $(\ref{H0_Robin})$, $(\ref{H1_Robin})$ and $(\ref{i-data})$, for $ N \geq 3$ and
\begin{displaymath}
p > 1 + \frac{2}{N},
\end{displaymath}
Problem $(\ref{Ext_Robin})$ admits global non-trivial positive solutions  for sufficiently small initial data $\varphi$.
\end{theo}
\emph{Proof: } Consider $\varphi$ satisfying $(\ref{i-data})$ and $v$ the non-trivial positive solution $v$ of the Neumann problem
\begin{equation*}
\left\{
    \begin{array}{ll}
        \partial _t v = \Delta v + v^p& \textrm{ in } \overline{\Omega} \times (0,+\infty), \\
        \partial _\nu  v= 0 & \textrm{ on }   \partial \Omega \times (0,+\infty), \\
        v(\cdot , 0) = \varphi & \textrm{ in } \overline{\Omega},
    \end{array}
\right.
\end{equation*}
where the initial data $\varphi$ is sufficiently small such that the solution $v$ is global. This choice can be achieved because $ N \geq 3$ and $p > 1 +  2/N$, see Levine \& Zhang \cite{LZ}. For all $\alpha \geq 0$ on $\partial \Omega \times(0,+\infty)$, we obtain
\begin{displaymath}
\partial _\nu v + \alpha v \geq 0 \textrm{ on } \partial \Omega \times (0,+\infty) .
\end{displaymath}
Thus, $v$ is a super-solution of Problem $(\ref{Ext_Robin})$, and we can deduce the statement of the theorem.
\\
\mbox{}\nolinebreak\hfill\rule{2mm}{2mm}\par\medbreak

\noindent Now, we suppose that there exists a positive constant $c>0$ such that
\begin{equation}\label{H2_Robin}
\alpha \geq c \textrm{ on } \partial \Omega \times \mathbb{R}^+ .
\end{equation}
We do not impose any condition on the dimension.

\begin{theo}
Let $\alpha$ be a coefficient satisfying $(\ref{H1_Robin})$ and $(\ref{H2_Robin})$, $\varphi$ an initial data with $(\ref{i-data})$. For
\begin{displaymath}
p > 1 + \frac{2}{N},
\end{displaymath}
Problem $(\ref{Ext_Robin})$ admits global positive solutions  for sufficiently small initial data $\varphi$.
\end{theo}
\emph{Proof: }
We consider the function $U :\overline{\Omega} \times [0,+\infty) \longrightarrow [0,+\infty)$ defined by
\begin{displaymath}
U(x,t) = A (t+t_0)^{-\mu} \exp \Big(  - \frac{\parallel x \parallel_2 ^2}{4(t+t_0)}\Big) ,
\end{displaymath}
where $\mu =1/(p-1)$, $t_0>0$  and $A>0$ will be chosen below. All the calculus will be detailed in the proof of the general Theorem \ref{R1_Robin}. If $A>0$ is small enough, we have
\begin{displaymath}
\partial_t U \geq \Delta U + U^p \textrm{ in } \overline{\Omega} \times (0,+\infty).
\end{displaymath}
On the boundary $\partial \Omega$, hypothesis $(\ref{H2_Robin})$ gives
\begin{eqnarray*}
\partial _\nu U(x,t) + \alpha U(x,t)  & \geq & \Big( \frac{ -x\cdot \nu(x)}{2(t+t_0)} + \alpha(x,t) \Big) U (x,t) \\
									& \geq & \Big( \frac{ -x\cdot \nu(x)}{2(t+t_0)} + c \Big) U(x,t) \\
\end{eqnarray*}
Since the boundary $\partial \Omega$ is compact, the function $(\partial \Omega \ni x \mapsto -x\cdot \nu(x) \in \mathbb{R}$ is bounded. We choose $t_0$ sufficiently big such that $ -x\cdot \nu(x)/(2t_0) + c \geq 0$. Then we obtain
\begin{displaymath}
\partial _\nu U + \alpha U \geq 0 \textrm{ on } \partial \Omega \times (0,+\infty).
\end{displaymath}
Finally, if we choose $\varphi \leq U(\cdot,0)$ in $\overline{\Omega}$, the function $U$ is a super-solution of Problem $(\ref{Ext_Robin})$.
\\
\mbox{}\nolinebreak\hfill\rule{2mm}{2mm}\par\medbreak

\begin{rema}
In the previous proof, one can note that the hypothesis $(\ref{H2_Robin})$ can be relaxed into 
\begin{equation}\label{Rem_Robin}
\alpha(x,t) \geq \frac{x\cdot \nu(x)}{2(t+t_0)} \textrm{ for all } (x,t) \in \partial \Omega \times (0,+\infty) .
\end{equation}
This condition gives us an optimal bound on $\alpha$ only if we know the geometry of the domain $\Omega$. For instance, if
\begin{displaymath}
\Omega = \{ \parallel x \parallel _2 > R \},
\end{displaymath}
we obtain
\begin{displaymath}
x\cdot \nu(x) = -R \textrm{ for all } x \in \partial \Omega.
\end{displaymath}
Then, the equation $(\ref{Rem_Robin})$ is equivalent to
\begin{equation*}
\alpha(x,t) \geq \frac{-R}{2(t+t_0)} \textrm{ for all } (x,t) \in \partial \Omega \times (0,+\infty) .
\end{equation*}
In particular, the previous theorem holds for all non-negative $\alpha$.
\end{rema}

\noindent In the one-dimensional case, using symmetry and translation, we can suppose that $\Omega = (-\infty,-1) \cup (1,+\infty)$. Then, without any additional hypothesis on the parameters of Problem $(\ref{Ext_Robin})$,  we obtain:

\begin{theo}
Assume the conditions $(\ref{H0_Robin})$, $(\ref{H1_Robin})$ and $(\ref{i-data})$. For dimension $N=1$ and
\begin{displaymath}
p > 3, 
\end{displaymath}
Problem $(\ref{Ext_Robin})$ admits global positive solutions  for sufficiently small initial data $\varphi$.
\end{theo}

\section{Generalization}

In the manner of Bandle \& Levine's results \cite{BL2}, we generalize our results. We consider the following problem
\begin{equation}\label{Ext_RobinG}
\left\{
    \begin{array}{ll}
        \partial _t u = L u + t^q \parallel x \parallel_2^s u^p& \textrm{ in }\overline{\Omega} \times (0,+\infty), \\
       \partial _\nu u+\alpha u= 0 & \textrm{ on } \partial \Omega \times (0,+\infty), \\
        u(\cdot , 0) = \varphi & \textrm{ in } \overline{\Omega} ,
    \end{array}
\right.
\end{equation}
where $q$ and $s$ are two positive real numbers, $p>1$ is a real number, and $L$ stands for the second order elliptic operator 
\begin{displaymath}
L= \sum_{i,j=1}^N \partial_{x_i} \Big( a_{ij}(x) \partial_{x_j} \Big ) +\sum_{i=1}^N b_i(x) \partial_{x_i} .
\end{displaymath}
To deal with classical solutions, the coefficients are assumed to be in $\mathcal{C}^2(\overline{\Omega})$. We keep the hypotheses $(\ref{H0_Robin})$, $(\ref{H1_Robin})$ and $(\ref{i-data})$  on the parameters $\alpha$ and $\varphi$. In order to state our principal results, we shall introduce some notations.
\begin{displaymath}
\rho(x)= \sum_{i,j=1}^N  a_{ij}(x) \frac{x_i x_j}{\parallel x \parallel_2^2}.
\end{displaymath}
Throughout, we assume that the matrix $A=(a_{ij})_{1\leq i,j\leq N}$ is normalized, so that for some $\nu_0 \in (0,1]$ 
\begin{displaymath}
0< \nu_0 \leq \rho \leq 1 \textrm{ in } \overline{\Omega}.
\end{displaymath}
Denote $b=(b_1, \dots , b_N)$ and let 
\begin{displaymath}
l(x)= \sum_{i,j=1}^N  \Big(\partial_{x_j} a_{ij}(x) - b_i(x)\Big)x_i ,
\end{displaymath}
\begin{displaymath}
l^*(x)= \sum_{i,j=1}^N  \Big(\partial_{x_j} a_{ij}(x) + b_i(x)\Big)x_i .
\end{displaymath}
We can state the following theorem concerning the blow-up case.

\begin{theo}
Assume that $N \geq 2$, 
\begin{displaymath}
\divergence b(x) \leq 0 \textrm{ in } \overline{\Omega} ,
\end{displaymath}
and
\begin{equation}\label{HG_Robin}
\rho(x) \leq \frac{\trace A(x) +l(x)}{2} \textrm{ in } \overline{\Omega} .
\end{equation}
Then, all non-trivial positive solutions of Problem $(\ref{Ext_RobinG})$ blow up in finite time for 
\begin{displaymath}
1 < p < 1 + \frac{2+2q+s}{N} .
\end{displaymath}
\end{theo}
\emph{Proof: }
Ab absurdo, we suppose that there exists a non-trivial positive solution $v$ of Problem $(\ref{Ext_RobinG})$. As in the proof of Theorem \ref{explode}, we deduce that there exists a non-trivial positive solution $u$ of the following Dirichlet problem
\begin{equation*}
\left\{
    \begin{array}{ll}
        \partial _t u = L u + t^q \parallel x \parallel_2 ^s u^p& \textrm{ in }\overline{\Omega} \times (0,+\infty), \\
        u= 0 & \textrm{ on } \partial \Omega \times (0,+\infty), \\
        u(\cdot , 0) = \varphi & \textrm{ in } \overline{\Omega} .
    \end{array}
\right.
\end{equation*}
According to Bandle \& Levine's results from \cite{BL2}, the solution $u$ blows up in finite time under the above hypotheses. Thus, $v$ must blow up too.
\\
\mbox{}\nolinebreak\hfill\rule{2mm}{2mm}\par\medbreak

\noindent For the one-dimensional case, Bandle \& Levine weaken the hypothesis $(\ref{HG_Robin})$. Then, we obtain:

\begin{theo}
Assume that $N = 1$, 
\begin{displaymath}
\divergence b(x) \leq 0 \textrm{ in } \overline{\Omega} ,
\end{displaymath}
and
\begin{displaymath}
\Big( \frac{2+2q+s}{p-1} -2 \Big) a_{11} +l >0 \textrm{ in } \overline{\Omega} .
\end{displaymath}
If $1 < p < 3+2q+s$, then all non-trivial positive solutions of Problem $(\ref{Ext_RobinG})$ blow up in finite time.
\end{theo}

\noindent Now, we consider the global existence case.

\begin{theo}\label{R1_Robin}
Assume that condition $(\ref{H2_Robin})$ is satisfied, 
\begin{displaymath}
\rho(x) \leq 1 \textrm{ in } \overline{\Omega} ,
\end{displaymath}
and
\begin{displaymath}
2\gamma_0 := \inf_{\overline{\Omega}} \Big( \trace A + l^* \Big) >0.
\end{displaymath}
Then, for any 
\begin{displaymath}
p>1+\frac{2+2q+s}{2\gamma_0} ,
\end{displaymath}
Problem $(\ref{Ext_RobinG})$ admits global non-trivial positive solutions if the initial data $\varphi$ is sufficiently small.
\end{theo}
\emph{Proof: }
We consider the function $U :\overline{\Omega} \times [0,+\infty) \longrightarrow [0,+\infty)$ defined by
\begin{displaymath}
U(x,t) = A (t+t_0)^{-\mu} \exp \Big(  - \frac{\parallel x \parallel_2 ^2}{4(t+t_0)}\Big) ,
\end{displaymath}
where $\mu =(2+2q+s)/(2p-2)$, $t_0>0$  and $A>0$ will be chosen below. We have
\begin{displaymath}
\partial_t U(x,t) = \Big( \frac{-\mu}{t+t_0}  + \frac{\parallel x \parallel_2 ^2}{4(t+t_0)^2}\Big)  U(x,t) ,
\end{displaymath}
\begin{displaymath}
L U(x,t) = \Big( \rho(x)\frac{\parallel x \parallel_2 ^2}{4(t+t_0)^2} -  \frac{\trace A + l^*}{2(t+t_0)} \Big)  U(x,t) ,
\end{displaymath}
and
\begin{displaymath}
\partial_\nu U(x,t) = \Big( \frac{-x\cdot\nu(x)}{2(t+t_0)}\Big)  U(x,t) .
\end{displaymath}
On the boundary $\partial \Omega$, we obtain:
\begin{displaymath}
\partial_\nu U(x,t) + \alpha U(x,t) = \Big( \frac{-x\cdot\nu(x)}{2(t+t_0)} +\alpha \Big)  U(x,t) .
\end{displaymath}
Thanks to hypothesis $(\ref{H2_Robin})$, and because the boundary $\partial \Omega$ is compact, we can choose $t_0$ sufficiently big such that
\begin{displaymath}
\frac{-x\cdot\nu(x)}{2t_0} +c  \geq 0 \textrm{ on } \partial \Omega.
\end{displaymath}
Thus, $\partial_\nu U(x,t) + \alpha U(x,t) \geq 0$ is achieved on $\partial \Omega \times (0,+\infty)$. Then, in $ \overline{\Omega}$, we have 
\begin{displaymath}
\partial_t U(x,t) - L U(x,t) = \Big( (1-\rho(x)) \frac{\parallel x \parallel_2 ^2}{4(t+t_0)^2} +  \frac{\trace A + l^* -2\mu }{2(t+t_0)} \Big)  U(x,t).
\end{displaymath}
With $\rho \leq 1$, we ignore the $t$-quadratic term, and by definition of $\gamma_0$, we obtain
\begin{equation}\label{Maj_Robin}
\partial_t U(x,t) - L U(x,t) \geq \Big(   \frac{ \gamma_0 -\mu }{t+t_0} \Big)  U(x,t),
\end{equation}
with $\gamma_0 -\mu>0$. On the other hand, $t<t+t_0$ implies
\begin{displaymath}
t^q \parallel x \parallel_2 ^s U^p(x,t) \leq  A^{p-1} \parallel x \parallel_2 ^s (t+t_0)^{q-\mu(p-1)}\exp\Big( \frac{-(p-1)\parallel x \parallel_2 ^2}{4(t+t_0)}\Big) U(x,t).
\end{displaymath}
Using the overestimation
\begin{displaymath}
\Big( \frac{2s}{p-1} \Big)^\frac{s}{2} \exp \Big( \frac{-s}{2}\Big) (t+t_0)^\frac{s}{2} \geq \parallel x \parallel_2^s \exp \Big( -\frac{\parallel x \parallel_2^2(p-1)}{4(t+t_0)}\Big) ,
\end{displaymath}
we obtain
\begin{equation}\label{Min_Robin}
t^q \parallel x \parallel_2 ^s U^p(x,t) \leq  A^{p-1} \Big( \frac{2s}{p-1} \Big)^\frac{s}{2} \exp \Big( \frac{-s}{2}\Big) (t+t_0)^{\frac{s}{2}+q-\mu(p-1)} U(x,t).
\end{equation}
By definition of $\mu$, we have $s/2+q-\mu(p-1)=-1$. Thus, we just have to choose $A$ sufficiently small, equations $(\ref{Maj_Robin})$ and $(\ref{Min_Robin})$ give
\begin{displaymath}
\partial_t U(x,t) - L U(x,t) \geq t^q \parallel x \parallel_2 ^s U^p(x,t).
\end{displaymath}
Finally, if the initial data $\varphi \leq U(\cdot,0)$ in $\overline{\Omega}$, $U$ is a super-solution of Problem $(\ref{Ext_RobinG})$, and we can deduce the existence of a solution using the truncation procedure of Section \ref{Prelem}.
\\
\mbox{}\nolinebreak\hfill\rule{2mm}{2mm}\par\medbreak


\begin{thebibliography}{10}

\bibitem{LZ} H.A. Levine and Q.S. Zhang, The critical Fujita number for a semilinear heat equation in exterior domains with homogeneous Neumann boundary values, Proc. of the R. Soc. of Edinb. 130A (200) 591-602.

\bibitem{Fujita} H. Fujita, On the blowing up of solutions of the Cauchy problem for $u_t =\Delta u + u^{1+\alpha}$, J. of the Fac. of Sci. of the Univ. of Tokyo 13 (1966) 109-124.

\bibitem{BL1} C.Bandle and H.A. Levine, On the existence and the nonexistence of global solutions of reaction-diffusion equations in sectorial domains, Trans. Am. Math. Soc. 316 (1989) 595-622.

\bibitem{BBR1} C. Bandle, J. von Below and W. Reichel, Parabolic problems with dynamical boundary conditions: eigenvalue expansions and blow up, Rend. Lincei Math. Appl. 17 (2006) 35-67.

\bibitem{JFR} J-F. Rault, The Fujita phenomenon in exterior domains under dynamical boundary conditions, Asymptot. Anal. 66 (2010) 1-8 .

\bibitem{BDC1} J. von Below and C. De Coster, A Qualitative Theory for Parabolic Problems under Dynamical Boundary Conditions, J. of Inequal. and Appl. 5 (2000) 467-486.

\bibitem{BL2} C. Bandle and H.A. Levine, Fujita type results for convective-like reaction diffusion equations in exterior domains, ZAMP 40 (1989) 665-676.

\bibitem{MS} K. Mochizuki  and R. Suzuki, Critical exponent and critical blow up for quasilinear parabolic equations, Isr. J. Math. 98 (1997) 141-156.

\bibitem{Suzuki} R. Suzuki, Critical blow-up for quasilinear parabolic equations in exterior domains, Tokyo J. of Math. 19 (1996) 397-409.

\end{thebibliography}
\end{document}